# On the history and use of some standard statistical models

## E. L. Lehmann[1]

*University of California, Berkeley*

**Abstract:** This paper tries to tell the story of the general linear model, which saw the light of day 200 years ago, and the assumptions underlying it. We distinguish three principal stages (ignoring earlier more isolated instances). The model was first proposed in the context of astronomical and geodesic observations, where the main source of variation was observational error. This was the main use of the model during the 19th century.

In the 1920's it was developed in a new direction by R.A. Fisher whose principal applications were in agriculture and biology. Finally, beginning in the 1930's and 40's it became an important tool for the social sciences. As new areas of applications were added, the assumptions underlying the model tended to become more questionable, and the resulting statistical techniques more prone to misuse.

## 1. Introduction

It was 200 years ago, in 1805, that Legendre first published the method of least squares and a vague formulation of what has come to be known as the standard linear model [27]. This model has played a central role in the statistical methodology used in the physical, biological and social sciences, and it is the aim of the present paper to sketch this role. I am not trying to write a history of the linear model[2], but am mainly concerned with the role the underlying assumptions have played in these three areas of application.

The model (defined in Section 6) assumes that each observation is the sum of two components: a deterministic term (which is a linear combination of the relevant explanatory variables) and a random term representing error or other "disturbances". The error terms are assumed to be independently, identically distributed (i.i.d.) according to a common normal distribution, and Sections 2–5 are therefore concerned with the role of these assumptions in the special case of the simple normal i.i.d. model. This model is of interest also in its own right as the standard model for the one-sample problem. Sections 6 and 7 take up issues concerning the non-random term of the linear model.

To conclude this introduction, let me briefly consider the general nature of mathematical, and particularly statistical, models.

---

[2]For accounts of this history see Seal [41] and Plackett [38], and in a broader context Stigler [45] and Hald [20].

The first use of the term "model" by a statistician that I have found occurs in Karl Pearson's Grammar of Science [33]. In this book he emphasizes the distinction between the real phenomenon and the model, which he usually calls "conception", but to which on a few occasions he refers as "models". For example he writes [33], p. 206:

> "The scientist postulates nothing in the world beyond sense [i.e. sense perceptions]; for him the atom and the ether are–like the geometric surface–models by which he resumes the world of sense.

The role and construction of mathematical models played a central role in the work of the mathematician Richard von Mises (1883-1953) who developed such models for a number of disciplines, among them aerodynamics, hydrodynamics, and plasticity. Besides undertaking these efforts in the physical sciences, von Mises felt that there was a great need for a similar treatment of probability theory. Instead of relying on Laplace's inadequate definition in terms of equally likely cases, he wanted to build a model for probability that would represent the physical reality underlying this concept.

In his fundamental paper of 1919 on the foundations of probability theory [48], he describes his approach as follows [my translation from the German]:

> "The present treatment is based on the assumption that probability theory is a *natural science, of the same kind as geometry or theoretical mechanics*. It has the aim to present the relations and dependencies of specific observable phenomenon, not as a faithful description of reality, but as its abstraction and idealization."

The observable phenomenon for which the paper builds a model is the stability of the long-run frequency of an outcome in a long sequence of repeated random events, for example, the frequency of heads in a long sequence of tosses of a coin. The principal constituents of the model are infinite sequences of trials with random outcomes, of which it is assumed that the frequency of a given outcome tends to a limit. To define randomness, von Mises requires that the same limit should obtain in any predetermined subsequence which can be chosen in the light of the observations up to this point.

It turned out that this formulation was too complicated; it was also too narrow in that it only applied to situations which allowed a large number of repetitions. As a formulation of probability theory it was replaced by a quite different approach proposed by Kolmogorov [23]. Kolmogorov did not construct a model for probability theory. Instead, he stated a small number of very simple axioms which any ideal model should satisfy. In this system probability itself was left undefined, subject only to these axioms. In this way the system could be interpreted and fleshed out not only by the frequency concept of probability but also by the idea of probability as degree of belief (The background of Kolmogorov's formulation is discussed in Shafer and Vovk [42]).

Although von Mises' modeling approach was not the one that was ultimately adopted by the profession, it nevertheless exerted great influence, particularly on Kolmogorov who cited the effect it had on his own formulation and who later took it as a starting point of a renewed effort to get a grip on the crucial (and difficult) concept of randomness[3].

---

[3]For a more detailed review of these foundational issues, see von Plato [49].

## 2. The assumption of normality in the 19th century

The most widely used statistical model for a sequence of repeated measurements has been the i.i.d. normal model according to which the observations are independently distributed with a common normal distribution. The first person to write down the formula for what today is called the normal density was De Moivre [10] who derived it as the limit of the binomial but did not consider it as a probability density in its own right. It was only in the beginning of the 19th Century through the combined insights of Laplace and Gauss that the central role of the normal distribution was realized. The normal distribution became the acknowledged model for the distribution of errors of physical (particularly astronomical) measurements and was called the Law of Errors.[4] It had a theoretical basis in the so called Law of Elementary Errors which assumed that an observational error is the sum of a large number of small independent errors and is therefore approximately normally distributed by the Central Limit Theorem. This argument was reinforced by extensive experience which showed good approximate agreement with the normal form.

The many textbooks on the subject all agreed on this point. For example, Brunt [6, 7], after citing theoretical arguments in favor of the normal distribution, states that "the final justification of Gauss' error curve rests upon the fact that it works well in practice and yields curves which in very many cases agree very closely with the observed frequency curves. The normal law is to be regarded as *proved by experience* and *explained* by Hagen's hypothesis [i.e. the law of elementary errors.]"

In a German text, Helmert [21, 22] explains [my translation]: "The form of the distribution of errors can only be determined through observation.... According to experience the [norma] law usually provides a close approximation."

Similar statements can be found in other British, German and French texts.

The history of the Theory of Errors is recounted in meticulous detail in a book on the subject by Czuber [9]. It was realized of course that the normal law could only be an approximation since in practice the observations were discrete and bounded. Nevertheless, as Czuber writes [my translation]:

> "An essential support of Gauss' Law of Errors is provided by the agreement which exists between its consequences and the results of observations that have really been obtained. It has led to the general acceptance by observers despite the concerns that can be raised against the various theoretical arguments in its favor."

Czuber devotes 14 pages to empirical comparisons of the law of errors with experience, and another section to the "elimination of contradictory observations" [i.e. gross errors], because it is only after such "doubtful observations" have been removed that the law applies.

He analyzes seven data sets from astronomy and geodesy, and in addition some experimental results of measurements carried out specifically to test the Gaussian law. His method in these examples is to compare the numbers of observations in various intervals with those predicted by theory. His conclusion is that the agreement is "satisfactory".

It should be noted that when Czuber mentions the general acceptance of the law, this statement is implicitly restricted to the principal areas of application he considers: the physical sciences and particularly astronomy and geodesy.

The normal distribution lost its exclusive position toward the end of the 19th Century when a strong interest developed in systematic application of statistical methods to biological, sociological and economic investigations. There the distributions encountered were often far from normal and frequently asymmetric. As a result, families of models were developed that include skewed and heavy-tailed distributions. The most influential of these at the time was the system proposed by Karl Pearson [34]. Stigler [45] (p. 335) explains Pearson's success:

> "In thirty pages of detailed examples he demonstrated the successful flexibility and practicality of the system with a force that bludgeoned any potential skeptic into submission. The examples ranged from Venn's barometric pressures to the heights of St. Louis school girls to Weldon's crabs to statistics on pauperism . . . . Not only was his work more general than others, it was practical and came to public view with a record of proven accomplishment."

---

[4]The term "normal distribution" was first suggested by Peirce [37], Lexis [28] and Galton [16] and began to take hold in the 1880's (see Stigler [46], Chapter 22).

## 3. Fisher and the assumption of normality

Early in the 20th Century the normal distribution regained its ascendancy through the small sample work of R.A. Fisher. As Geary [17] summarized the situation in a historical overview:

> "Our historian will find a significant change of attitude a quarter of a century ago following the brilliant work of R.A. Fisher who showed that when universal normality could be assumed, inferences of the widest practical usefulness could be drawn. Prejudice in favor of normality returned in full force and interest in non-normality receded in the background."

But could universal normality be assumed? Fisher developed his new methods not in the context of astronomical measurements but for application to biological and agricultural data, and here the assumption is on much shakier ground.

In the earlier applications, the principal source of variation in the observations was observational error. Now to this is added the variability of the subjects (trees, agricultural plots, farms, . . . ) being measured. These subjects Fisher modeled as random samples from an infinite population, and he based his derivations on the assumption that the population distribution of the characteristic being measured is normal. But this distribution depends on the nature of the population, and the basis for such an assumption in these circumstances often is weak.

Fisher's treatment of the assumptions of normality in his enormously successful 1925 book, "Statistical Methods for Research Workers" (SMRW) [11], gave rise to a heated controversy between Fisher and E.S. Pearson. The conflict had its origin in Pearson's 1929 review of the 2nd Edition of SMRW for the journal "Nature". On the whole the review was favorable, but it contained the following critical paragraph.

> "There is one criticism however which must be made from the statistical point of view. A large number of these tests are based . . . on the assumption that the population sampled is of the 'normal' form. That this is the case can be gathered from a careful reading of the test, but the point is not sufficiently emphasized. It does not appear reasonable to lay stress on the 'exactness' of the tests when no means whatever are given of appreciating how rapidly they became inexact as the population sampled diverges from normality. That the tests, for example, connected with the analysis of variance [here he is referring to the $F$-test for variances] are far more dependent on normality than those involved Student's $z$ (or $t$) distribution is almost certain, but no clear indication of the need for caution in their application is given."

To make things worse, the review also contained the sentence: "It would seem wiser in the long-run even in a textbook, to admit the incompleteness of theory in this direction, rather than risk giving the reader the impression that the solution of all his problems has been achieved."

In defense of his criticism,[5] Pearson (in a letter to Gosset) cited a paper by the American economist Tolley [47] which showed that its author indeed had been misled by Fisher's book. Tolley wrote:

> "Recently the English School of Statisticians has developed formulas and probability tables to accompany them which, they state, are applicable regardless of the form of the frequency distribution. These formulas are given, most of them without proof, in Fisher's book [11] .... If we accept the statements of those who have developed these newer formulas, skew frequency distributions and small samples need cause no further difficulty as far as measurement of error is concerned."

Fisher was shaken by Tolley's misunderstanding and in a letter to Gosset (June 27, 1929) admitted some culpability:

> "The claim of exactness for the solutions and tests given was wrong, although a careful reader would find that I had kept within the letter of the law by hidden allusions to normality."

A careful reading of SMWR concerning its treatment of normality in fact shows the following.

1. For some of the results, the assumption is stated (though never emphasized).
2. More often, results are stated without any qualifications.
3. The treatment of examples is very deficient in this respect. As Gosset pointed out to Fisher (June 24, 1929):

   > "Although when you think about it you agree that 'exactness' or even appropriate use depends on normality, in practice you don't consider the question at all when you apply your tables to your examples: *not one word.*"

4. In the chapter on distributions (which introduces the normal, Poisson and binomial and no others), Fisher states that it is important to know "the experimental conditions upon which they occur." However, neither this nor any later chapter contains any discussion of the conditions under which data can be expected to be normally distributed (not one word, as Gosset might say.)
5. Another important omission (not even taken up in later editions after it had been pointed out by Pearson) is the great sensitivity to the assumption of normality in tests of variances (rather than means), which make these tests essentially useless.

This criticism must of course be seen against the background of the enormous achievements and novelty of the book. In the process of making accessible Fisher's research of the preceding decade, it established a new paradigm and revolutionized statistical methodology.

The book provides little information about Fisher's own attitude toward the assumption of normality, which forms the basis of his work on the analysis of variance, covariance and regression. The clearest statement of his position that I have been able to find is in a letter of 1929 to "Nature" entitled "Statistics and biological research", which was the concluding document in the controversy with Pearson. Basing his defense of the assumption on experience rather than theory, Fisher claims:

> "On the practical side there is little enough room for anxiety, especially among biologists, who are used to checking the adequacy of their methods by control experiments. The difficulty of obtaining decisive results often flows from heterogeneity of materials, often from causes of bias, often, too, from the difficulty of setting up an experiment in such a way as to obtain a valid estimate of errors. I have never known difficulty to

---

[5]For more details on this controversy see Pearson [32] and the Fisher-Gosset correspondence (Gosset [18]).

> arise in biological work from imperfect normality of the variation, often though I have examined data for this particular cause of difficulty; nor is there, I believe, any case to contrary in the literature. This is not to say that the deviation from 'Student's' *t*-distribution found by Shewhart and Winters [43], for samples from rectangular and triangular distributions, may not have a real application in some technological work, but rather that such deviations have not been found, and are scarcely to be looked for, in biological research and ordinarily conducted."

Fisher is not so naive as to claim universal normality even for the kind of biological data with which he is dealing. What he claims instead is that his methods are insensitive to departures from normality, in modern terminology that they are fairly robust against non-normality. This turned out, as Pearson had already discovered through simulation studies, to be true for tests of means but not for tests of variances. (However, even for means a theorem of Bahadur and Savage [3] suggest the need for some caution. They showed that for any given sample size, no matter how large, there exist distributions for which the size of the *t*-test is arbitrarily close to 1.)

## 4. The role of normality after Fisher

Fisher's "Statistical Methods" was enormously influential. Together with the work on which it was based it was the most important instrument for the movement from 19th Century large-sample to 20th Century small-sample statistics. The volume sold well, particularly considering how small the statistical community was at the time. A second edition became necessary in 1928 and a third in 1930, with the first three editions selling 1050, 1250 and 1500 copies respectively. This flow continued throughout Fisher's life with new editions appearing every 2 to 3 years.

However, in the 1930's other texts also began to bring the new methodology to ever widening audiences. By far the most successful of these was George Snedecor's "Statistical Methods" [44], which was published in 1937, and the seven editions of which sold the unparalleled number of 237,000 copies. For several years it was one of the most cited publications in the Science Citation Index, and in 1995 still had nearly 2000 entries.[6]

Snedecor's book is essentially a more clearly written, simpler, very user-friendly version of Fisher's SMRW. It also contains a large number of numerical examples, mainly from agriculture and biology. The book avoids complications and pays practically no attention to the assumptions underlying the recommended procedures including the assumption of normality. This is in line with Snedecor's philosophy which he explains in the Preface:

> "To the mathematical statistician must be delegated the task of developing the theory and devising the methods, *accompanying these latter by adequate statements of the limitations of their use* (my italics)."

He adds that "None but the biologist can decide whether the conditions are fulfilled in his experiments." This seems begging the question because how can the biologists or other users make this decision when they are not clearly told what are the conditions.

This lack of attention to assumptions is shared by most of the many texts on statistics in the 40's and 50's that followed those of Fisher and Snedecor. A fairly comprehensive look at the role of the assumptions underlying Fisher's small-sample methods was provided in 1959 by Scheffé in his book "The Analysis of Variance"

---

[6]Cited from Carriquiri and David [8].

[40], in a justly famous chapter of nearly 40 pages entitled "The effects of departures from the underlying assumptions."

Scheffé measures the departure of a distribution from normality in terms of the coefficients $\gamma_1$ of skewness and $\gamma_2$ of kurtosis, the standardized 3rd and 4th moments. For any symmetric distribution $\gamma_1 = 0$, and for the normal distribution also $\gamma_2 = 0$. In the introduction to the chapter, he provides an idea of the range of $(\gamma_1, \gamma_2)$ in engineering data, in routine chemical analyses, in ages of marriage, in barometric heights, and in the length and breadth of beans.

Scheffé points out that for large $n$ the distribution of Student's $t$ is approximately normal and hence "the inferences about the mean which are valid in the case of normality must be correct for large $n$ regardless of the form of the population." This is followed by a discussion of the effect of non-normality on the $\chi^2$-test for variance, where it is shown that non-normality causes serious errors. Thus, Scheffé concludes that

> "the effect of violation of the normality assumption is slight on inferences about the mean but dangerous on inferences about variances."

He points out that these results had already been noted by E. S. Pearson [31], and their reason by Box [4].

However, Scheffé's concern is an exception. The standard textbook attitude toward the assumption of normality is well summarized by Brownlee [5] (p. 179):

> "It is presumably on this foundation [that means are approximately normal by the Central Limit Theorem] that applied statisticians have found empirically that usually there is no great need to fuss about the normality assumption. After a statistician has analyzed several quite widely differing transformations of a variable in a fair number of specific instances and found that the conclusions reached are substantially identical for all the transformations, then he ceases to worry unduly about the normality assumption in most situations."

## 5. The assumption of independence

The standard model for the one-sample problem, i.e. for a number of repeated measurements or other observations of a common quantity, assumed not only that the observations are normally distributed but also that they are independent. This assumption has received much less attention than the assumption of normality. As Kruskal [26] states: "An almost universal assumption in statistical models for repeated measurements of real-word quantities is that these measurements are independent, yet we know that such independence is fragile."

As an explanation of the casual assumption of independence he suggests:

> "One answer is ignorance.... Far more important than simple ignorance, is seductive simplicity: It is so easy to multiply marginal probabilities, formulas simplify, and manipulation is relatively smooth, so the investigator neglects dependence, or hopes that it makes little difference. Sometimes the hope is realized, but more often dependence can make a tremendous difference."

In many situations, the assumption of independence seems natural, even obvious, because of the absence of any direct influence of one observation on another. However, a "spurious" (but nevertheless very real) dependence may be caused by the presence of a common factor. This phenomenon was investigated (and the term "spurious correlation" coined) by Karl Pearson in two papers of 1897 and 1902 [35, 36], in which he reports extensive experiments on measurements carried out by himself.

A difficulty for any serious investigation of the affect of dependence is the great variety of forms that dependence can take. Even Scheffé, who includes a discussion of dependence in his chapter on failures of assumptions (mentioned in the preceding section) cannot do more than examine two or three special cases. His conclusion: "The effect of correlation in the observations can be very serious on inferences about means."

Following Scheffé, many of the later books on linear models and regression analysis included discussion of the effects of non-normality, dependence, and inequality of variances. Particularly noteworthy in this regard is Miller [30] who provides a careful treatment of these issues as an integral part of each chapter of his book.

There was also a trickle-down effect on some of the more general introductions to statistical methods. Thus, the seventh (1980) edition of Snedecor's text (Snedecor had died in 1974 and Cochran had become a co-author) contained a substantial chapter on "Failure in the assumptions". Its treatment of dependence closely paralleled that of Scheffé and in particular also warned about the effect of dependence on the *t*-test.

## 6. The linear model

The previous sections have been concerned with the i.i.d. normal model for the one-sample problem, i.e. for repeated measurements or observations of the same quantity. In the remainder of the paper we shall consider the general linear model where the same issues arise as well as some new ones.

The model is given by

$$Y_i = \sum_{j=1}^n x_{ij}\beta_j + \varepsilon_i \quad (i = 1, \dots, n)$$

where the $Y$'s are the observed values, the $x$'s are known constants, the $\beta$'s unknown parameters, and the $\varepsilon$'s the errors. Of the $\varepsilon$'s it is typically assumed that they are independently normally distributed with mean 0 and common variance $\sigma^2$.

This model and the proposal to estimate the $\beta$'s by means of least squares, are due to Gauss and Legendre at the beginning of the 19th Century. Known as the theory of combinations of observations, the resulting methodology, applied primarily to astronomy and geodesy, became the principal statistical activity throughout much of the 19th Century, with many textbooks devoted to it. (When Gosset in 1904 needed statistical methods for his brewery work, his main sources of information were Airy's [1] (3rd Ed) book with the ponderous title "On the Algebraic and Numerical Theory of Errors of Observations" and Merriman's "A Textbook on Least Squares" [29]. As Stigler [45] (p.11) puts it: "The method of least squares was the dominant theme – the leitmotif – of nineteenth-century mathematical statistics."

In this theory, probability calculations such as that of determining the probable error of the estimates were carried out assuming large samples. In this way the variance of the normal error distribution could be assumed known, and the estimates (which were linear functions of the observations) could therefore be assumed to be (approximately) normally distributed with known variance.

The linear model became a much more flexible instrument through Fisher's introduction of analysis of variance and covariance, and regression analysis. At the same time he extended its application to the more complicated data of biology with their complex sources of variation mentioned in Section 4. Fisher was able to overcome many of the difficulties arising in this extension through the use of randomization and other aspects of experimental design set forth in his 1935 book on this subject [12].

After Fisher's biological applications, the linear model still faced one important challenge: the application of regression models to the social sciences, particularly to economics. Here in addition to the assumptions concerning the random terms $\varepsilon_i$ (normality and independence), questions concerning the adequacy of the structural part $\Sigma x_{ij}\beta_j$ became particularly important.

Freedman [15] (p. 9) compares the situation with that in astronomy where

> "The relevant variables were known from Newtonian mechanics, and so were the functional forms of the equations connecting them. Measurement could be done with high precision. Much was known about the errors in the measurements and the equations."

Regarding the corresponding issues in the social sciences, I shall here restrict attention to economics[7] where after some early isolated instances, statistical model building was pursued systematically (beginning in the 1930's) by writers such as Tinbergen, Haavelmo [19] and Koopmans. It was developed further in the 1940's by the work of the Cowles Commission, which was then the center of economic research. About these efforts Koopmans [24] points out that

> "this theory has been widely applied to data obtained from agricultural experiments or from measurements in biological populations. There are some essential differences between data of this kind and those usually encountered in economic problems.
>
> In agricultural experiments some of the determining variables can be completely controlled by the experimenter . . . . Other determining variables less under his control are usually by their nature subject to adequate independent variation. In that respect they bear a resemblance to the variables representing measurable characteristics of individuals of a biological population which is usually conceived as random drawings from a *stable* (my italics) probability distribution.
>
> In economic analysis variables at the control of an experimenting institution are exceptional. Further only a few types of variables . . . are so erratic in nature that they could reasonably be regarded as drawings from any stable distribution . . . .
>
> Further the relations between the variables studied in this type of analysis are themselves subject to gradual or abrupt changes, according to institutional or technical changes in society . . . . Therefore the number of observations from which the regression coefficients have to be estimated is limited by the very nature of the problem."

Commenting on the relative roles of the economist and the statistician in this work, Koopmans states:

> "The economist – or, in general, the expert in the field to which the dependent variable belongs – should by economic reasoning and general economic experience . . . devise a set of determining variables which he expects to be a complete set [i.e. such that the effect of any additional variables can safely be absorbed into the error term]."

Koopmans makes it clear that the determination of the appropriate structural part is not the task of the statistician but of the subject matter expert. This point is also made by Arrow [2] when he states that

> "The method of scientific investigation indicated in the preceding paragraphs calls then for intensive a priori thinking to formulate a model, followed by the selection of a best-fitting structure from that model by appropriate statistical techniques."

Once the model has become formulated, statistical techniques make it possible to deduce far-reaching conclusions. However the reliability of these conclusions is of course limited by the reliability of the model. The difficulty of verifying the model assumptions on the one hand, and the power of the statistical machinery to deliver results on the other, have provided a fertile ground for misuse of the regression methodology. The chief critic of this misuse has been David Freedman who over the past 25 years has called attention to it, both in a series of close to 50 substantial papers, and as a statistical consultant and expert witness in more than 100 cases.

---

[7]For applications to psychology see Krüger et al. [25] (Chapters 2 and 3).

As an example of Freedman's criticism consider his detailed analysis of a book by Hope, which deals with the effects of education on class mobility (Freedman [13]). In this paper Freedman points out that

> "one problem noticeable to a statistician is that investigators do not pay attention to the stochastic assumptions behind the models. It does not seem possible to derive these assumptions from current theory, nor are they easily validated empirically on a case-by-case basis"

The paper ends with the devastating conclusion:

> "My opinion is that investigators need to think more about the underlying process, and look more closely at the data, without the distorting prism of conventional (and largely irrelevant) stochastic models. Estimating nonexistent parameters cannot be very fruitful. And it must be equally a waste of time to test theories on the basis of statistical hypotheses that are rooted neither in prior theory nor in fact, even if the algorithms are recited in every statistics text without caveat."

The opening sentence of this conclusion suggests an alternative to the methodology Freedman is criticizing: better data, more substantive knowledge and input, multiple studies under varying conditions as had been carried out for example in establishing smoking as a major cause of lung cancer and other diseases. These requirements had been foreshadowed in the passages from Koopmans and Arrow cited earlier in this section.

Freedman capped his long involvements in these issues as scientist, consultant and expert witness with a text: "Statistical Models–Theory and Practice" [15]. As the title suggests, the book provides an account not only of the theory of statistical modeling but also–through the careful examination of many real-life examples–a guide to how such modeling should and should not be used.

An additional valuable feature of the book is a review of the literature on statistical modeling which is both a resource and encouragement for further study.

## 7. Conclusions

This paper has been concerned with the assumptions underlying the linear model and the special case of the normal i.i.d. model. It considered three assumptions: normality, independence and in the more general case, the linear structure of the deterministic part. Of these three, the assumption of normality has received by far the most attention in the literature. For the i.i.d. case and a few simple linear models it has led to an alternative nonparametric methodology which has been developed to avoid it. Unfortunately these nonparametric methods are no help with respect to dependence (or the inequality of variances, a topic we have not discussed here.)

A point emphasized by this survey is that different fields of study involve different kinds of data and have very different modeling situations. Thus, when the only source of variation is observational error, as often is the case in the physical sciences, the situation is much simpler than in the biological sciences where one is dealing with samples from a somewhat heterogeneous (but stable) population. And these situations, in turn, are easier to handle than those of observational studies in the social sciences where the populations are less stable and where random sampling typically is not possible.

This subject matter dependence creates a difficulty for general texts on statistical methods but it does not absolve them from responsibility. As a minimum, such texts should provide warnings against the unsubstantiated use of standard models. A useful cautionary example might be the $F$-test of variances with its strong dependence on the assumption of normality. Even with the $t$-test which is fairly robust against non-normality one should not tempt students into Tolley's error (discussed in Section 4).

A warning that should be included in such texts is the danger of the too facile assumption of independence. Another distinction that would be worth mentioning is that between observational and experimental studies, perhaps with references to the literature (for example the book by Rosenbaum [39]) treating the special methods developed for the former.

Though general statistics texts could, and should include, this kind of material, they can only go so far. As mentioned by Koopmans, the statistician and the subject matter specialist each have a role to play. A general text cannot provide the subject matter knowledge and the special features that are needed for successful modeling in specific cases. Experience with similar data is required, knowledge of theory and, as Freedman points out: shoe leather.

**Acknowledgments.**   I am grateful to Persi Diaconis and Juliet Shaffer for many helpful critical comments.

## References

[1]  AIRY, G. (1861, 1879). *On the Algebraical and Numerical Theory of Errors of Observations*, 3rd ed. Macmillan, London.

[2]  ARROW, K. (1951). Mathematical models in the social sciences. In *The Policy Sciences* (Lerner and Lasswell, eds.). Stanford Univ. Press. MR0044815

[3]  BAHADUR, R. AND SAVAGE, L. (1956). The nonexistence of certain statistical procedures in nonparametric problems. *Ann. Math. Statist.* **27** 1115–1122. MR0084241

[4]  BOX, G. (1953). Non-normality and tests of variances. *Biometrika* **40** 318–335. MR0058937

[5]  BROWNLEE, K. (1960). *Statistical Theory and Methodology in Science and Engineering.* Wiley, New York. MR0119268

[6]  BRUNT, D. (1917). *The Combination of Observations.* Cambridge Univ. Press.

[7]  BRUNT, D. (1931). *The Combination of Observations.* Cambridge Univ. Press.

[8]  CARRIQUIRY, A. AND DAVID, H. (2001). George Waddel Snedecor. In *Statisticians of the Centuries* (Heyde and Seneta, eds.). Springer, New York.

[9]  CZUBER, E. (1891). *Theorie der Beobachtungsfehler.* Teubner, Leipzig.

[10]  DE MOIVRE, A. (1733). Approximatio ad Summam Terminorum Binomii $(a+b)^n$ in Seriem Expansi. Printed for private circulation.

[11]  FISHER, R. (1925). *Statistical Methods for Research Workers.* Oliver and Boyd, Edinburgh.

[12]  FISHER, R. (1935). *The Design of Experiments.* Oliver and Boyd, Edinburgh.

[13]  FREEDMAN, D. (1987). As others see us (with discussion). *J. Educ. Statist.* **12** 101–223. Reprinted in J. Shaffer, ed. The Role of Models in Nonexperimental Social Science AERA/ASA Washington, D.C. (1997).

[14]  FREEDMAN, D. (1991). Statistical models and shoe leather. In *Sociological Methodology* (P. Marsden, ed.). Amer. Social Assoc., Washington, D.C.

[15] FREEDMAN, D. (2005). *Statistical Models: Theory and Practice.* Cambridge Univ. Press.

[16] GALTON, F. (1877). Typical laws of heredity. *Nature* **15** 492–495, 512–514, 532–533.

[17] GEARY, R. (1947). Testing for normality. *Biometrika* **34** 209–242. MR0023497

[18] GOSSET, W. (1970). Letters from W. S. Gosset to R. A. Fisher, 1915–1936 with summaries by R. A. Fisher and a Foreword by L. McMullen. Printed for private circulation.

[19] HAAVELMO, T. (1944). The probability approach in econometrics. *Econometrica* **12** Supplement. MR0010953

[20] HALD, A. (1998). *A History of Mathematical Statistics from 1750 to 1930.* Wiley, New York. MR1619032

[21] HELMERT, F. (1872). *Die Ausgleichsrechnung nach der Methode der Kleinsten Quadrate.* Teubner, Leipzig.

[22] HELMERT, F. (1907). *Die Ausgleichsrechnung nach der Methode der Kleinsten Quadrate.* Teubner, Leipzig.

[23] KOLMOGOROV, A. (1933). *Grundbegriffe der Wahrscheinlichkeitsrechnung.* Springer, Berlin.

[24] KOOPMANS, T. (1937). *Linear Regression Analysis of Economic Time Series.* Netherlands Economic Institute, Haarlem.

[25] KRÜGER, L., GIGERENZER, G. AND MORGAN, M., EDS. (1987). *The Probabilistic Revolution* **2**. MIT Press, Cambridge, MA. MR0929869

[26] KRUSKAL, W. (1988). Miracles and statistics: The casual assumption of independence. *J. Amer. Statist. Assoc.* **83** 929–940. MR0997572

[27] LEGENDRE, A. (1805). *Nouvelles Méthodes pour la Détermination des Orbites des Comètes.* Courcier, Paris.

[28] LEXIS, W. (1877). *Theorie der Massenerscheinungen in der Menschlichen Gesellschaft.* Wagner, Freiburg.

[29] MERRIMAN, M. (1884). A textbook on the method of least squares, 8th ed. 1900.

[30] MILLER, R. (1986). *Beyond ANOVA, Basics of Applied Statistics.* Wiley, New York. MR0838087

[31] PEARSON, E. (1931). The analysis of variance in cases of non-normal variation. *Biometrika* **23** 114–133.

[32] PEARSON, E. (1990). *Student.* Clarendon Press, Oxford. MR1255103

[33] PEARSON, K. (1892). *Grammar of Science.* Walter Scott, London. 3rd ed. of 1911 reprinted by Meridian Books, 1957.

[34] PEARSON, K. (1895). Contributions to the mathematical theory of evolution, II. Skew variation in homogeneous material. *Phil. Trans. Roy. Soc. London* **186** 343–414.

[35] PEARSON, K. (1897). On a form of spurious correlation which may arise when indices are used in the measurement of organs. *Proc. Roy. Soc.* **60** 489–497.

[36] PEARSON, K. (1902). On the mathematical theory of errors of judgement, with special reference to the personal equation. *Phil. Trans. Roy. Soc. London* **198** 235–299.

[37] PEIRCE, C. (1873). On the theory of errors of observations. Appendix 21 to the Report of the Superintendent of the U.S. Coast Survey for the year ending June 1870. 200–224. Reprinted in *Stigler: American Contributions to Mathematical Statistics in the Nineteenth Century* **2** (1980). Arno Press, New York.

[38] PLACKETT, R. (1972). The discovery of the method of least squares. *Biometrika* **59** 239–251. MR0326871

[39] ROSENBAUM, P. (1995, 2002). *Observational Studies.* Springer, New York.

[40] SCHEFFÉ, H. (1959). *The Analysis of Variance.* Wiley, New York. MR0116429

[41] SEAL, H. (1967). The historical developmnent of the gauss linear model. *Biometrika* **54** 1–24. MR0214170

[42] SHAFER, G. AND VOVK, V. (2006). The sources of kolmogorovŌs *Grundbe-griffe. Statist. Sci.* **21** 70–98. MR2275967

[43] SHEWHART, W. AND WINTERS, F. (1928). Small samples – new experimental results. *J. Amer. Statist. Assoc.* **23** 144–153.

[44] SNEDECOR, G. (1937). *Statistial Methods.* The Iowa State College Press, Ames, Iowa.

[45] STIGLER, S. (1986). *The History of Statistics.* Belknap Press, Cambridge, MA. MR0852410

[46] STIGLER, S. (1999). *Statistics on the Table.* Harvard Univ. Press, Cambridge, MA. MR1712969

[47] TOLLEY, H. (1929). Economic data from the sampling point of view. *J. Amer. Statist. Assoc.* **24** 69–72.

[48] VON MISES, R. (1919). Grundlagen der Wahrscheinlichkeitsrechnung. *Math. Zeitschrift* **5** 52–99.

[49] VON PLATO, J. (1994). *Creating Modern Probability.* Cambridge Univ. Press.